%% file: width_bound-article2.tex
\documentclass{article}

\usepackage[leqno,fleqn]{amsmath}
\usepackage{amsthm}
\usepackage{amsbsy}
\usepackage{comment}
\usepackage{enumerate}
\usepackage[normalem]{ulem}
\usepackage[titletoc,title]{appendix}

\usepackage{tikz}
\usetikzlibrary{positioning,patterns,calc}

\usepackage{amsfonts}
\usepackage{mathrsfs}
\usepackage{MnSymbol}
\DeclareMathAlphabet{\mathpzc}{OT1}{pzc}{m}{it}

\newcommand{\duline}[1]{{\bgroup \markoverwith{{\bgroup \markoverwith{\rule[-1.2pt]{0.1pt}{0.4pt}}\ULon {\rule[-2.8pt]{1pt}{0.4pt}}}}\ULon {#1}}} 
\newcommand{\dulineF}[1]{{\bgroup \markoverwith{{\bgroup \markoverwith{\rule[-1.2pt]{0.4pt}{0.4pt}}\ULon {\rule[-2.8pt]{2pt}{0.4pt}}}}\ULon {#1}}} 
\newcommand{\dulineeta}[1]{{\bgroup \markoverwith{{\bgroup \markoverwith{\rule[0.2pt]{0.1pt}{0.4pt}}\ULon {\rule[-1.4pt]{1pt}{0.4pt}}}}\ULon {#1}}} 

\newcommand{\real}{\mathbb{R}}
\newcommand{\integer}{\mathbb{Z}}
\newcommand{\naturals}{\mathbb{N}}
\newcommand{\supp}{\textup{supp}}
\newcommand{\haus}{\mathcal{H}}
\newcommand{\rest}[1]{\,\rule[-.23cm]{.012cm}{.5cm}_{\,#1}} 
\newcommand{\compactf}[1]{C_c(#1)} 
\newcommand{\grassp}[2]{G(#1,#2)} 
\newcommand{\grassm}[2]{G_{#1}(#2)} 
\newcommand{\restset}{\,\!\raisebox{-0.4ex}{\ensuremath{\invbackneg}}\,} 
\newcommand{\bdry}[1]{\textup{bdry}\left(#1\right)} 
\newcommand{\interior}[1]{\textup{int}\left(#1\right)} 
\newcommand{\closure}[1]{\overline{#1}} 

\renewcommand{\span}{\textup{span}}
\renewcommand{\index}{\textup{index}}
\newcommand{\nullity}{\textup{null}}
\newcommand{\vol}{\textup{vol}}

\newcommand{\unitvol}{\alpha} 
\newcommand{\mass}{\textbf{\duline{M}}} 
\newcommand{\fmetric}{\textbf{\dulineF{F}}} 
\newcommand{\fnorm}{\mathcal{F}} 
\newcommand{\rectv}{\uuline{\upsilon}} 
\newcommand{\tcone}{\mathscr{C}} 
\newcommand{\varf}{\mathcal{V}} 
\newcommand{\rect}{\textup{R}\mathcal{V}} 
\newcommand{\intv}{\textup{I}\mathcal{V}} 
\newcommand{\replacement}{\mathcal{R}} 
\newcommand{\acurr}{\mathfrak{A}} 
\newcommand{\cycles}{\mathcal{Z}} 
\newcommand{\rectc}{\uuline{\mathfrak{t}}} 
\newcommand{\width}{\omega} 
\newcommand{\sweepout}{\mathcal{P}} 
\newcommand{\domain}{\textup{dmn}}

\theoremstyle{plain}
\newtheorem{theorem}{Theorem}[section]
\newtheorem{prop}[theorem]{Proposition}
\newtheorem{lemma}[theorem]{Lemma}
\newtheorem{corollary}[theorem]{Corollary}
\newtheorem{claim}{\texttt{\textbf{\underline{Claim}}}}
\makeatletter
\@addtoreset{claim}{theorem}
\makeatother

\theoremstyle{definition}
\newtheorem{definition}[theorem]{Definition}

\theoremstyle{remark}
\newtheorem*{remark}{\underline{Remark}}
\numberwithin{equation}{section}

\begin{document}

\title{The width of Ellipsoids}
\author{Nicolau Sarquis Aiex}
\renewcommand{\footnotemark}{}
\thanks{The author was supported by a CNPq-Brasil Scholarship}
\maketitle

\renewcommand{\abstractname}{\vspace{-\baselineskip}}
\begin{abstract}
\noindent \textsc{Abstract.} We compute the `$k$-width of a round $2$-sphere for $k=1,\ldots,8$ and we use this result to show that unstable embedded closed geodesics can arise with multiplicity as a min-max critical varifold.
\end{abstract}

\input{wb-intro-article}

\hfill

\textit{
Acknowledgements: I am thankful to my PhD adviser Andr\'e Neves for his guidance and suggestion to work on this problem.
}

\input{wb-comp-article} 

\input{wb-network-article}

\input{wb-almostm-article2}

\input{wb-ellip-article}

\input{wb-problems-article}

\input{wb-appendix-article}


\input{wb-bib}
\end{document}

%% file: wb-intro-article.tex
\section{Introduction}
The aim of this work is to compute some of the $k$-width of the $2$-sphere.
Even in this simple case the full width spectrum is not very well known.
One of the motivations is to prove a Weyl type law for the width as it was proposed in \cite{mgromov1}, where the author suggests that the width should be considered as a non-linear spectrum analogue to the spectrum of the laplacian.

In a closed Riemannian manifold $M$ of dimension $n$ the Weyl law says that $\frac{\lambda_p}{p^{\frac{2}{n}}}\rightarrow c_n\vol(M)^{-\frac{2}{n}}$ for a known constant $c_n$, where $\lambda_p$ denotes the $p^{th}$ eigenvalue of the laplacian.
In the case of curves in a $2$-dimensional manifold $M$ we expect that
\begin{equation*}
   \frac{\width_p}{p^{\frac{1}{2}}}\rightarrow C_2\vol(M)^{-\frac{1}{2}},
\end{equation*}
where $C_2>0$ is some constant to be determined.

We were unable to compute all of the width of $S^2$ but we propose a general formula that is consistent with our results and the desired Weyl law.
In the last section we explain it in more details.

By making a contrast with classical Morse theory one could ask the following two naive questions about the index and nullity of a varifold that achieves the width:
\hfill

\noindent \textbf{\underline{Question 1:}} Let $(M,g)$ be a Riemannian manifold and $V\in\intv_l(M)$ be a critical varifold for the $k$-width $\width_k(M,g)$. Then
\begin{equation*}
   k\leq \index(V)+\nullity(V).
\end{equation*}
\noindent \textbf{\underline{Question 2:}} Let $(M,g)$ be a Riemannian manifold and $V\in\intv_l(M)$ be a critical varifold for the $k$-width $\width_k(M,g)$. Then
\begin{equation*}
   \index(V)\leq k.
\end{equation*}
Where $\index(V)$ and $\nullity(V)$ are the index and nullity of the second variation $\delta^{(2)}V$ on the space of vectorfields in $M$.
By a critical varifold we mean that $V$ is obtained as the accumulation point of a minmax sequence.

As a pertubation of our results we will show that \textbf{\underline{Question 1}} is false for one-varifolds on a surface.
Regarding \textbf{\underline{Question 2}}, it was recently shown by Marques-Neves in \cite{fmarques-aneves3} that $\index(\supp(V))\leq k$ in the case of codimension one and $3\leq dim(M)\leq 7$.
The authors also conjecture that the two-sided unstable components of $V$ must have multiplicity one.
In the hypersurface case the Pitts' min-max theorem gives us an embedded minimal hypersurface, whereas the dimension $1$ case allows self-intersections.
That is why they do not expected it to hold for curves.
This work provides a concrete example of how it fails to be true in the dimension $1$ case.

To illustrate these questions we present an example in which it holds and explain how it fails in our context.
Say we are trying to study closed geodesics by analyzing the energy functional $E$ in the free loop space $\Lambda=W^{1,2}(S^1,M)$, in which case we can apply infinite dimensional Morse theory.
Take $a<b$ regular values and suppose we can find a non-trivial homology class $\alpha\in H_k(\Lambda^b,\Lambda^b)$ ($\Lambda^a=\{E\leq a\}$) then we can find a closed geodesic $\gamma$ satisfying
\begin{equation*}
E(\gamma)=\inf_{A\in\alpha}\sup_{x\in \supp\,A}E(x).
\end{equation*}
In this case it is known that $\index(\gamma)\leq k\leq \index(\gamma)+\nullity(\gamma)$ (this is encoded in \cite[\textsection 1 Lemma 2]{dgromoll_wmeyer1}, alternatively see \cite[Chapter 2 Corollary 1.3]{kchang1}).
Compared to our case $\gamma$ would correspond to $V$, a non-trivial $k$-dimensional homology class corresponds to a $k$-sweepout and the minmax quantity is analogue to the $k$-width.

There are two differences between the classical Morse Theory set up and Almgren-Pitts minmax.
The first is that we are working with varifolds instead of parametrized curves, which allow degenerations.
On the other hand we compute the index and nullity in the same way, by using vectorfield variations.

As an example, consider the union of two great circles in the $2$-sphere.
It divides the sphere into four discs and for each of them we take a $1$-parameter contraction to a point.
If we follow the boundary of these contractions simultaneously we would have a $1$-parameter family of cycles that decreases length.
However, this is not generated by an ambient vectorfield, so it does not contribute to the index of the stationary varifold.

The other difference is that Almgren-Pitts minmax theory works with homotopies instead of homologies, which forces us to consider different variations to obtain the critical varifold.


This article is divided as follows.
In section 2 we briefly overview definitions and main properties of sweepouts, currents and varifolds. 
In section 3 we define geodesic networks, that will be the candidates of critical varifolds for the width.
Here we prove a structure result for $1$-dimensional stationary integral varifolds.
In section 4 we define almost minimising varifolds and characterize the singularities of such varifolds.
This will allow us to have a regularity result for critical varifolds of low parameter widths.
In section 5 we compute the $k$-width of $S^2$ for $k=1,\ldots,8$.
Then we use the regularity results to find the critical varifolds for a generic ellipsoid.
Though we could not explicitly show which width correspond to each critical varifold we prove that it provides a counterexample anyway.

%% file: wb-comp-article.tex
\section{Preliminaries}

In this section definitions and notations are established. 
Throughout this section $M$ denotes a closed Riemannian manifold of dimension $m$ isometrically embedded in $\real^n$ for some $n>0$.

Let us denote by $\cycles_k(M)$ the space of \textbf{flat $k$-cycles} in $M$ with coefficients in $\integer_2$ endowed with the flat topology.
We write $\fnorm$ for the \textbf{flat norm} and $\mass$ for the \textbf{mass} of a cycle. 

We adopt the definition of varifolds in \cite{jpitts1}.
We denote the spaces of $k$-\textbf{varifolds}, \textbf{rectifiable varifolds} and \textbf{integral varifolds} by $\varf_k(M)$, $\rect_k(M)$ and $\intv_k(M)$, respectively.
These spaces are endowed with the weak topology induced by the metric $\fmetric$.

Given a rectifiable varifold $V\in\rect_k(M)$ we write $\tcone_pV$ for the tangent cone of $V$ at the point $p\in\supp\|V\|$.
We also denote by $\grassp{k}{n}$ the space of $k$-planes in $\real^n$ and $\grassm{k}{M}=\{(x,P)\in\real^n\times\grassp{k}{n}:x\in M, P\subset T_xM\}$ the $k$-Grassmanian bundle over $M$.
For a rectifiable set $S\subset \real^n$ and $\theta$ and integrable function in $\grassm{k}{\real^n}$ we write $\rectv(S,\theta)$ the varifold associated to $S$ with density $\theta$.

Now we establish a relation between currents and varifolds.
Given a $k$-current $T$ (not necessarily closed) we denoted by $|T|\in\varf_k(M)$ the varifold induced by the support of $T$ and its coefficients.
Reversely, given a $k$-varifold $V$ we denote by $[V]$ the unique $k$-current such that $\Theta^k(|[V]|,x)=\Theta^k(V,x)\,mod\,2$ for all $x\in\supp\|V\|$ (see \cite{bwhite2}).

\subsection{Sweepouts and the width}

In \cite{falmgren1} Almgren proved, in particular, that $\pi_i(\cycles_k(M))=H_{i+k}(M;\integer_2)$ for all $i>0$.
We call it the \textbf{Almgren isomorphism} and denote it by $F_A$.
It follows from the Universal Coefficient Theorem that $H^{n-k}(\cycles_k(M);\integer_2)=\integer_2$, denote its generator by $\bar\lambda$ and $\bar\lambda^p$ the cup product with itself $p$ times.
For the next definition we follow \cite{lguth1} and \cite{fmarques-aneves1}.
\begin{definition}
Let $X\subset [0,1]^N$ be a cubical subcomplex for some $N>0$ and $f:X\rightarrow\cycles_k(M)$ a flat continuous map.
We say that $f$ is a $p$-\textbf{sweepout} if
\begin{equation*}
   f^*(\bar\lambda^p)\neq 0 \in H^{p(n-k)}(X;\integer_2).
\end{equation*}
Denote the set of $p$-sweepouts with \textit{no concentration of mass} (see definition \ref{no mass concentration}) in $M$ by $\sweepout_p(M)$.

We define the $p$-\textbf{width} of $(M,g)$ as
\begin{equation*}
   \width_p(M,g)=\inf_{f\in\sweepout_p}\sup_{x\in\domain(f)}\mass(f(x)),
\end{equation*}
where $\domain(f)$ denotes the domain of $f$.
\end{definition}

Note that $\width_p\leq\width_{p+1}$ since every $(p+1)$-sweepout is also a $p$-sweepout.

\subsection{Varifolds in $S^n$}
Let $(S^n,g_{S^n})$ denote the round sphere of radius $1$ in $\real^{n+1}$. Given a varifold $V\in\varf_k(S^n)$ we can define the cone generated by $V$ in $\real^{n+1}$. It is sufficient to define a positive functional in the space $\compactf{\grassm{k+1}{\real^{n+1}}}$ (see for example \cite[\textsection 5.2(3)]{wallard1}).

\begin{definition}
Given $V\in\varf_k(S^n)$ define $C(V)\in\varf_{k+1}(\real^{n+1})$ to be the measure corresponding to the functional
\begin{equation*}
   C(V)(f)=\int_0^{\infty}\tau^{k}V(f_{\tau})d\tau
\end{equation*}
where $f\in \compactf{\grassm{k+1}{\real^{n+1}}}$ and $f_{\tau}\in \compactf{\grassm{k}{\real^{n+1}}}$ is given by
\begin{equation*}
   f_{\tau}(x, P) = \left\{
   \begin{aligned}
      & f\left(\tau x, P\oplus\real\langle x\rangle \right) , \textup{ if } x\in S^n \textup{ and } P\subset T_xS^n; \\
      & 0 , \textup{ otherwise.}
   \end{aligned}
\right.
\end{equation*}
\end{definition}

%
%
\begin{prop}\label{coneproperties}
The cone map $C:\varf_{k}(S^n)\rightarrow\varf_{k+1}(\real^{k+1})$ satisfy the following properties:
\begin{enumerate}[(i)]
   \item $C(V)$ is a cone varifold; \label{conepropertiescone}
   \item If $a,b\in\real_{\geq 0}$ and $V,W\in\varf_k(S^n)$ then $C(aV+bW)=aC(V)+bC(W)$; \label{conepropertieslinear}
   \item If $V\in\rect_k(S^n)$ then $C(V)$ is given by
   \begin{equation*}
      C(V)(f)=\int_0^{\infty}\tau^k\int_{S^n} f\left(\tau x, T_{V}(x)\oplus\real\langle x\rangle\right)\Theta^k(V,x)d\haus^k_x d\tau,
   \end{equation*}
   where $f\in\compactf{\grassm{k+1}{\real^{n+1}}}$ and $T_V(x)\subset T_x S^n$ is the tangent space of $V$ defined $\|V\|$-almost everywhere in $S^n$; \label{conepropertiesformula}
\end{enumerate}
\end{prop}

\begin{proof}
\textbf{(\ref{conepropertiescone}):} We must show that ${\eta_{0,\lambda}}_\# C(V)=C(V)$ for all $\lambda>0$, where $\eta_{0,\lambda}(x)=\lambda x$ for $x\in\real^{n+1}$. 
Take any $f\in\compactf{\grassm{k+1}{\real^{n+1}}}$ and compute
\begin{equation*}
   \begin{aligned}
      {\eta_{0,\lambda}}_\# C(V)(f) & = \int f(x,\tilde P) d({\eta_{0,\lambda}}_\# C(V))_{x,\tilde P}\\
                                    & = \int [J^{k+1}\eta_{0,\lambda}](x)f(\eta_{0,\lambda}(x),D\eta_{0,\lambda}\rest{x}\cdot \tilde P)dC(V)_{x,\tilde P}\\
                                    & = \frac{1}{\lambda^{k+1}}\int f\left(\frac{x}{\lambda},\tilde P\right)dC(V)_{x,\tilde P}\\
                                    & = \frac{1}{\lambda^{k+1}}\int_0^\infty\tau^{k}V(f_{\frac{\tau}{\lambda}})d\tau\\
                                    & = \frac{1}{\lambda^{k+1}}\int_0^\infty (t\lambda)^{k}V(f_t)\lambda dt\\
                                    & = C(V)(f).
   \end{aligned}
\end{equation*}
Here it was used the definition of pushforward and change of variables $t=\frac{\tau}{\lambda}$ in the second last line.\\
\textbf{(\ref{conepropertieslinear}):} This is straightforward from the definition.\\
\textbf{(\ref{conepropertiesformula}):} To prove this formula simply use that
\begin{equation*}
V(f_\tau)=\int_{S^n}f_\tau\left(x,T_V(x)\right))d\|V\|_x
\end{equation*}
holds for rectifiable varifolds and $d\|V\|_x=\Theta^k(V,x)d\haus^k_x$.
\end{proof}

\indent We now want to prove that this cone map is continuous with respect to the weak convergence. 

%
%
\begin{lemma}\label{conecontinuity}
Let $\{V_n\}\subset\varf_k(S^n)$ be a sequence of varifolds converging to $V\in\varf_k(S^n)$ in the $\fmetric$-metric. Then $C(V_n)\rightarrow C(V)$ with respect to $\fmetric$.
\end{lemma}

\begin{proof}
It is enough to prove that $C(V_n)(f)$ $\rightarrow C(V)(f)$ for any compactly supported function in $\grassm{k+1}{\real^{n+1}}$.

There exist $R_0>0$ such that $\supp(f)\subset B(0,R_0)\times\grassp{k+1}{n+1}$. For $\tau>R_0$ we have
\begin{equation*}
   f_\tau(\tau x,P)=0,
\end{equation*}
for all $x\in S^n$ and $P\in\grassp{k}{n+1}$. Thus, whenever $\tau>R_0$,
\begin{equation*}
   V_n(f_\tau)=\int f(\tau x, P)d\left(V_n\right)_{x, P}=0
\end{equation*}
for all $n>0$, and $V(f_\tau)=0$. This implies that the sequence $h_n(\tau)=\tau^k V_n(f_\tau)$ is uniformly bounded. By the Dominated Convergence theorem we obtain
\begin{equation*}
   \begin{aligned}
      \lim_{n\rightarrow\infty}C(V_n)(f) & = \lim_{n\rightarrow\infty}\int_0^{R_0} h_n(\tau)d\tau=\int_0^{R_0}\lim_{n\rightarrow\infty}h_n(\tau)d\tau\\
                                         & = \int_0^{R_0} \tau^kV(f_\tau)d\tau\\
                                         & = C(V)(f).
   \end{aligned}
\end{equation*}
\end{proof}

Next we show that the cone of a varifold associated to a rectifiable set in $S^n$ is defined by the cone of the set, as one would expect.
%
%
\begin{lemma}\label{conerectifiable}
Let $R\subset S^n$ be a $k$-rectifiable set, $\theta:\grassm{k}{S^n}\rightarrow\real_{\geq 0}$ a locally integrable function and $\rectv(R,\theta)\in\rect_k(S^n)$. Then
\begin{equation*}
   C\left(\rectv(R,\theta)\right)=\rectv\left(\tilde{R}, \tilde{\theta}\right),
\end{equation*} 
where $\tilde{R}=\{\lambda x\in\real^{n+1}|\lambda\geq 0,\; x\in R\}$ and $\tilde{\theta}:\grassm{k+1}{\real^{n+1}}\rightarrow\real_{\geq 0}$ is a given by
\begin{equation*}
   \tilde\theta(x,\tilde{P})=\left\{
   \begin{aligned}
      & \theta\left(\frac{x}{|x|},P\right),\; \textup{if }x\neq 0 \textup{ and } \tilde{P}=P\oplus\langle x\rangle;\\
      & 0,\; \textup{otherwise}.
   \end{aligned}
   \right.
\end{equation*}
\end{lemma}

\begin{proof}
It is easy to see that $\tilde\theta$ is locally integrable in $\grassm{k+1}{\real^{n+1}}$, $\tilde R$ is $(k+1)$-rectifiable and its tangent space is given by $T_{\tilde R}(x)=T_R(\frac{x}{|x|})\oplus\langle\frac{x}{|x|}\rangle$ for $x$ $(\haus^{k+1}\restset\tilde{R})$-almost-everywhere. 
For $f\in\compactf{\grassm{k+1}{\real^{n+1}}}$ compute
\begin{equation}\label{compt1}
   \begin{aligned}
      \rectv\left(\tilde{R},\tilde{\theta}\right)(f) & = \int_{\real^{n+1}}f\left(x,T_{\tilde{R}}(x)\right)\tilde{\theta}(x,T_{\tilde{R}}(x))d(\haus^{k+1}\restset\tilde{R})_x\\
                                                     & = \int_{\real^{n+1}\setminus\{0\}}f\left(x,T_R\left(\frac{x}{|x|}\right)\oplus\langle \frac{x}{|x|}\rangle\right)\theta\left(\frac{x}{|x|},T_R\left(\frac{x}{|x|}\right)\right)d(\haus^{k+1}\restset\tilde{R})_x.
   \end{aligned}\tag{*}
\end{equation}

\indent We want to use the Co-area formula (see \cite[\textsection 3.2.22]{hfederer1}), we clarify notation and make some remarks.\\
\indent Define the warped product metric on $(0,+\infty)\times S^n$ as $g_{(\tau,x)}=d\tau^2+\tau^2(g_{S^n})_x$, where $g_{S^n}$ is the round Riemannian metric on $S^n$. Let $d_g$, $d_{S^n}$ and $d_{0}$ be the metrics induced by $g$ on $(0,+\infty)\times S^n$, $g_{S^n}$ on $S^n$ and the Euclidian metric $g_0$ on $\real^{n+1}$ respectively. Given any metric $d$ we denote by $\haus^k(d)$ the $k$-dimensional Hausdorff measure associated to $d$.

\begin{claim} The metrics $g$, $g_{S^n}$ and $g_0$ satisfy:
\begin{enumerate}[(a)]
   \item $F:(\tau,z)\in\left((0,+\infty)\times S^n, g\right)\mapsto \tau z\in(\real^{n+1}\setminus \{0\}, g_0)$ is an isometry;
   \item $d_g\left((\tau,z),(\tau,y)\right)=\tau d_{S^n}(z,y)$;
   \item $(\iota_{\tau})_*\haus^k(d_{S^n})=\tau^{-k}\haus^k(d_g)$, where $\iota_\tau:S^n\rightarrow (0,+\infty)\times S^n$ is the inclusion in the slice $\{\tau\}\times S^n$.
\end{enumerate}
\end{claim}
Firstly, (a) is a well known fact and (b) follows easily from the definition. Lastly, (b) implies that $\iota_\tau$ is $\tau^{-1}$-Lipschitz so (c) follows from basic properties of $\haus^k$.

For simplicity denote $h(x)$ the integrand in (\ref{compt1}). Applying a change of variables and the Co-area formula for the projection $(\tau,z)\mapsto \tau$ we obtain
\begin{equation*}
   \begin{aligned}
      \rectv\left(\tilde{R},\tilde{\theta}\right)(f) & = \int_{\real^{n+1}\setminus\{0\}}h(x)dF_*F^{-1}_*(\haus^{k+1}(d_0)\restset\tilde{R})_x\\
                                                     & = \int_{(0,+\infty)\times S^n}h\circ F(\lambda,z)d\left(\haus^{k+1}(d_g)\restset(0,+\infty)\times R\right)_{(\lambda,z)}\\
                                                     & = \int_{(0,+\infty)\times R}h\circ F(\lambda,z)d\haus^{k+1}(d_g)_{(\lambda,z)}\\
                                                     & = \int_0^{\infty}\left(\int_{\{\tau\}\times R}h\circ F(\lambda,z)d\haus^k(d_g)_{(\lambda,z)}\right)d\haus^1(d_0)_\tau.
   \end{aligned}
\end{equation*}
Changing variables again and using (c) we conclude
\begin{equation*}
   \begin{aligned}
      \rectv\left(\tilde{R},\tilde{\theta}\right)(f) & = \int_0^{\infty}\left(\int_{\{\tau\}\times R}h\circ F(\lambda,z)d\left(\tau^k{\iota_{\tau}}_*\haus^k(d_{S^n})\right)_{(\lambda,z)}\right)d\haus^1(d_0)_\tau\\
                                                     & = \int_0^\infty\left(\int_R h\circ F\circ \iota_\tau (z) \tau^k d\haus^k(d_{S^n})_z\right)d\tau\\
                                                     & = \int_0^\infty\tau^k\int_{S^n} h(\tau z)d(\haus^k\restset R)_z d\tau.
\end{aligned}
\end{equation*}
The proof is finished by replacing $h$ in the formula and \ref{coneproperties}(iii).
\end{proof}

Finally, we can prove the main properties of the cone $C(V)$.

%
%
\begin{prop}\label{conemain}
Let $V\in\rect_k(S^n)$ and $C(V)\in\varf_{k+1}(\real^{n+1})$. Then the following is true:
\begin{enumerate}[(i)]
   \item $C(V)$ is rectifiable;\label{conemainrectifiable}
   \item if $V$ is integral then so is $C(V)$;\label{conemainintegral}
   \item $\supp\|C(V)\|=\{\lambda x\in\real^{n+1}|x\in\supp\|V\| \; and \; \lambda\geq 0\}$;\label{conemainsupport}
   \item if $y\neq 0$ then $\Theta^{k+1}(C(V),y)=\Theta^k(V,\frac{y}{|y|})$;\label{conemaindensity}
   \item if $y\neq 0$ then
      \begin{equation*}
         \|V\|(S^n)=(k+1)\lim_{r\rightarrow\infty}\frac{\|C(V)\|B(y,r)}{r^{k+1}};
      \end{equation*}\label{conemainmass}
   \item if $V$ is stationary and $k\geq 1$ then so is $C(V)$.\label{conemainstationary}
\end{enumerate}
\end{prop}

\begin{proof}
\textbf{(\ref{conemainrectifiable}):} Let $V = \lim_{n\rightarrow\infty}\sum_{i=1}^n\rectv(R_i,\theta_i)$, where $R_i\subset S^n$ is $k$-rectifiable and $\theta_i:\grassm{k}{S^n}\rightarrow \real_{\geq 0}$ is locally integrable for all $i>0$. 
The result follows directly from\ref{coneproperties}(\ref{conepropertieslinear}), Lemma \ref{conecontinuity} and Lemma \ref{conerectifiable}.

\textbf{(\ref{conemainintegral}):} Just note in the proof of Lemma \ref{conerectifiable}, if $\theta$ is integer-valued then so is $\tilde\theta$.

\textbf{(\ref{conemainsupport}):} First show that $\supp \|C(V)\|\supset \{\lambda x\in\real^{n+1}|x\in\supp\|V\| \; and \; \lambda\geq 0\}$. Take $y\not\in\supp\|C(V)\|$ and a positive continuous function $\tilde f:\real^{n+1}\rightarrow\real_{\geq 0}$ supported in $B(y,r)$, for some $r>0$. Define $f(x)=\tilde{f}(a x)$, where $a=\min\{1,|y|\}$. So $f$ is supported in $B(\frac{y}{|y|},r)$. If we assume $\|C(V)\|(\tilde f)=0$ then it is easy to check that $\|V\|(f)=0$, so $\frac{y}{|y|}\not\in \supp\|V\|$. The other inclusion is similar.

\textbf{(\ref{conemaindensity}):} for simplicity put $C=C(V)$. It is enough to show that
\begin{equation*}
\int_{\real^{n+1}\setminus\{0\}} g(y)d\|C\|_y=\int_{\real^{n+1}\setminus\{0\}} g(y)\Theta^k\left(V,\frac{y}{|y|}\right) d\haus^{k+1}_y
\end{equation*}
for every continuous function $g$ compactly supported in $\real^{n+1}\setminus\{0\}$.

\indent If $f\in\compactf{\grassm{k+1}{\real^{n+1}}}$ satifies $f(0,\tilde P)=0$ for all $\tilde P\in\grassp{k+1}{n+1}$ then, from rectifiability (property (\ref{conemainrectifiable})), it follows that
\begin{equation*}
   \begin{aligned}
   C(f) & = \int_{\grassm{k+1}{\real^{n+1}}}f(y,\tilde P)dC_{(y,\tilde P)}\\
        & = \int_{\real^{n+1}\setminus\{0\}} f(y,T_C (y))d\|C\|_y.
   \end{aligned}
\end{equation*}
On the other hand, by property \ref{coneproperties}(\ref{conepropertiesformula}), we have
\begin{equation*}
   C(f)  = \int_0^\infty\tau^k\int_{S^n} f(\tau x, T_V(x)\oplus\real\langle x\rangle)\Theta^k(V,x)d\haus^k_xd\tau.
\end{equation*}
Following a computation similar to the proof of Lemma \ref{conerectifiable} we conclude
\begin{equation*}
   C(f) = \int_{\real^{n+1}\setminus\{0\}}f\left(y,T_V\left(\frac{y}{|y|}\right)\right)\Theta^k\left(V,\frac{y}{|y|}\right)d\haus^{k+1}_y
\end{equation*}
\indent Take a continuous function $g:\real^{n+1}\rightarrow \real$ compactly supported in $\real^{n+1}\setminus\{0\}$ and define $f(y,\tilde P)=g(y)$ for all $y\in\real^{n+1}$ and $\tilde P\in\grassp{k+1}{n+1}$.
The result follows by replacing such $f$ in the previous formulas.

\textbf{(\ref{conemainmass}):} Fix $y\neq 0$, $r>0$ and let $F:\real^{n+1}\setminus\{0\}\rightarrow (0,+\infty)\times S^n$ be the isometry $F(z)=(|z|,\frac{z}{|z|})$.
Denote $A(r)=F(B(y,r))$, $\textup{pr}_1(\tau,x)=\tau$.

If $r>|y|$ then $0\in B(y,r)$ and
\begin{equation*}
   \begin{aligned}
   \|C\|B(y,r) & = \int_{B(y,r)}d\|C\|_z\\
               & = \int_{B(y,r)\setminus\{0\}}\Theta^{k+1}(C,z)d\haus^{k+1}(d_0)_z\\
               & = \int_{B(y,r)\setminus\{0\}}\Theta^{k}\left(V,\frac{z}{|z|}\right)d\haus^{k+1}(d_0)_z\\
               & = \int_{A(r)}\Theta^{k}(V,x)d\haus^{k+1}(d_g)_{(\tau,x)}.
   \end{aligned}
\end{equation*}
Furthermore, $\textup{pr}_1(A(r))=(a(r),b(r))$, with $a(r)=\inf_{z\in B(y,r)}|z|=0$ and $b(r)=\sup_{z\in B(y,r)}|z|=|y|+r$. 
Note also that $\textup{pr}_1^{-1}(\tau)=\{\tau\}\times S^n$ for $\tau< r-|y|$.
Applying the Co-area formula with respect to $\textup{pr}_1$we get
\begin{equation*}
   \begin{aligned}
      \|C\|B(y,r) & = \int_{0}^{|y|+r}\int_{\textup{pr}_1^{-1}(\tau)}\Theta^k(V,x)d\haus^k(d_g)_{(\lambda,x)}d\tau\\
                  & = \int_{0}^{r-|y|}\int_{\{\tau\}\times S^n}\Theta^k(V,x)d\haus^k(d_g)_{(\lambda,x)}d\tau \\
                  & \quad + \int_{r-|y|}^{r+|y|}\int_{\textup{pr}_1^{-1}(\tau)}\Theta^k(V,x)d\haus^k(d_g)_{(\lambda,x)}d\tau\\
                  & =  \int_{0}^{r-|y|}\tau^k\int_{S^n}\Theta^k(V,x)d\haus^k(d_{S^n})_{(x)}d\tau \\
                  & \quad + \int_{r-|y|}^{r+|y|}\int_{\textup{pr}_1^{-1}(\tau)}\Theta^k(V,x)d\haus^k(d_g)_{(\lambda,x)}d\tau.
   \end{aligned}
\end{equation*}
The first term in the sum is given by
\begin{equation*}
\int_{0}^{r-|y|}\tau^k\int_{S^n}\Theta^k(V,x)d\haus^k(d_{S^n})_{(x)}d\tau = \frac{(r-|y|)^{k+1}}{k+1}\|V\|(S^n).
\end{equation*}
Since $\textup{pr}_1^{-1}(\tau)\subset \{\tau\}\times S^n$, the second term is bounded by
\begin{equation*}
   \int_{r-|y|}^{r+|y|}\int_{\textup{pr}_1^{-1}(\tau)}\Theta^k(V,x)d\haus^k(d_g)_{(\lambda,x)}d\tau \leq \frac{(r+|y|)^{k+1}-(r-|y|)^{k+1}}{k+1}\|V\|(S^n).
\end{equation*}
When we divide by $r^{k+1}$ and take the limit $r\rightarrow\infty$, the first term converges to $\frac{\|V\|(S^n)}{k+1}$ and the second term tends to zero.

\textbf{(\ref{conemainstationary}):} Since we assume $k \geq 1$ it is enough to prove that $C(V)$ is stationary outside the origin.

Fix a vector field $Y$ with compact support $\supp (Y)\subset \real^{n+1}\setminus \{0\}$. We can write $Y(y)=h(y)y+X(\frac{y}{|y|})$ where $X$ is a compactly supported vector field in $S^n$ and $h$ is a compactly supported function. The first variation is given by $\delta C(Y)=\delta C (h(y)y)+\delta C (X(\frac{y}{|y|}))$. Let us compute the first term:
\begin{equation*}
   \begin{aligned}
   \delta C(h(y)y) & = \int_{\grassm{k+1}{\real^{n+1}}}\textup{div}_{\tilde P}(h(y)y)dC_{(y,\tilde P)}\\
                   & = \int_{0}^{\infty}\tau^k\int_{S^n}\textup{div}_{T_V(x)\oplus\real\langle x \rangle}(h(\tau x) \tau x)d\|V\|_xd\tau\\
                   & = \int_{0}^{\infty}\tau^k\int_{S^n}Dh\rest{\tau x}\cdot\tau x + h(\tau x)\textup{div}_{T_V(x)\oplus\real\langle x \rangle}(\tau x) d\|V\|_xd\tau\\
                   & = \int_{0}^{\infty}\int_{S^n}\tau^k\left(\tau Dh\rest{\tau x}\cdot x+h(\tau x)(k+1)\right)d\|V\|_xd\tau\\
                   & = \int_{S^n}\left[\int_0^{\infty}\left(\frac{d}{dt}\rest{t=\tau}t^kh(t x)\right)d\tau \right]d\|V\|_x\\
                   & = 0
   \end{aligned}
\end{equation*}
In the last line we used that $h$ has compact support away from $0$.\\ 
\indent Using that $X$ doesn't depend on the radial direction, that is, $\textup{div}_{\langle x\rangle}(X)=0$, we compute the second term
\begin{equation*}
   \begin{aligned}
   \delta C(X(\frac{y}{|y|})) & = \int_{\grassm{k+1}{\real^{n+1}}}\textup{div}_{\tilde P}(X(\frac{y}{|y|}))dC_{(y,\tilde P)}\\
                              & = \int_{0}^{\infty}\tau^k\int_{S^n}\textup{div}_{T_V(x)\oplus\real\langle x \rangle}(X(x))d\|V\|_xd\tau\\
                              & = \int_{0}^{\infty}\tau^k\int_{S^n}\textup{div}_{T_V(x)}(X)+\textup{div}_{\real\langle x\rangle}(X)d\|V\|_x d\tau\\
                              & = \int_{0}^{\infty}\tau^k\delta V(X)d\tau\\
                              & = 0.
   \end{aligned} 
\end{equation*}
Thus finishing the proof of the proposition.
\end{proof}

%% file: wb-network-article.tex
\section{Geodesic Networks}
In this section we are concerned with $1$-dimensional varifolds whose support is represented by geodesic segments. 
Our aim is to prove that any stationary integral $1$-varifold has this structure.

\begin{definition}
Let $U\subset M$ be an open subset. 
A varifold $V\in\intv_1(M)$ is called a \textbf{geodesic network} in $U$ if there exist geodesic segments $\{\alpha_1,\ldots,\alpha_l\}$ in $M$ and $\{\theta_1,\ldots,\theta_l\}\subset \integer_{>0}$ such that
\begin{enumerate}[(a)]
   \item
      \begin{equation*}
         V\restset\grassm{k}{U} = \sum_{j=1}^l\rectv(\alpha_j\cap U,\theta_j).
      \end{equation*}
   \item Let $\Sigma_V=\cup_{j=1}^{l}(\partial\alpha_j)\cap U$, we require each $p\in\Sigma_V$ to belong to exactly $m=m(p)\geq 3$ geodesic segments $\{\alpha_{j_1},\ldots\alpha_{j_m}\}$ and
\begin{equation*}
   \sum_{k=1}^m\theta_{j_k}\dot{\alpha}_{j_k}(0)=0.
\end{equation*}
   Here we are taking the arc-length parametrization with start point at $p$.
\end{enumerate}
We call a point in $\Sigma_V$ a \textbf{junction}. 
We say that a junction is \textbf{singular} if there exist at least $2$ geodesic segments with $\theta_{j_k}\dot{\alpha_{j_k}}(0)\neq-\theta_{j_{k'}}\dot{\alpha_{j_{k'}}}(0)$ and \textbf{regular} otherwise.
A \textbf{triple junction} is a point $p\in\Sigma_V$ such that $p$ is the boundary of only $3$ geodesic segments with multiplicity $1$ each.
\end{definition}

%
%
%
%
%
\indent The following properties can be derived straightforwardly from the definition.
%
%
\begin{prop}\label{networkproperties}
Let $V$ be a geodesic network in $U\subset M$. The following holds:
\begin{enumerate}[(i)]
   \item $V$ is stationary in $U$;\label{networkstationary}
   \item if $p\in\Sigma_V$ and $\{(\alpha_{j_1},\theta_{j_1}),\ldots,(\alpha_{j_m},\theta_{j_m})\}$ define this junction then the tangent cone at $p$ is given by
   \begin{equation*}
      \tcone_pV=\sum_{k=1}^m\rectv(\textup{cone}(\dot{\alpha}_{j_k}(0)),\theta_{j_k})
   \end{equation*}
   where $\textup{cone}(\dot{\alpha}_{j_k}(0))=\{\lambda\dot{\alpha}_{j_k}(0)\in T_pM| \lambda\geq 0\}$ and $\alpha_{j_k}(0)=p$.\label{networkcone}
\end{enumerate}
\end{prop}

%
%
\begin{corollary}\label{networktriple}
Let $U\subset M$ be an open set. If $V$ is a geodesic network in $U$ and $\Theta^1(V,x)<2$ for all $x\in\supp\|V\|$, then every $p\in\Sigma_V$ is a triple junction.
\end{corollary}
\begin{proof}
First note that the condition $\Theta^1(V,x)<2$ at regular points imply that $\theta_j=1$ for all $j$. By proposition \ref{networkproperties}(\ref{networkcone}) the density is given by
\begin{equation*}
   \Theta^1(V,p)=\Theta^1(\tcone_pV,0)=\sum_{k=1}^m\frac{\theta_{j_k}}{2}.
\end{equation*}
Since $\theta_{j_k}= 1$ we must have $m < 4$ thus $m=3$.
\end{proof}

In the two dimensional case we can infer further on the regularity of junctions.

\begin{corollary}\label{characterization}
Let $M$ be a surface and $V\in\intv_1(M)$ be a geodesic network with density $\Theta^1(V,p)\leq 2$ for all $p\in\supp\|V\|$. Then \textbf{either}
\begin{enumerate}[(i)]
   \item $\Sigma_V$ contains at least one triple junctions \textbf{or}
   \item $\Sigma_V$ has no triple junctions, all junctions are regular and $V$ is given by 
      \begin{equation*}
         V=\sum_{i=1}^l\rectv(\gamma_i,1)
      \end{equation*}
      where $\gamma_i$ are closed geodesics (possibly repeated) and $\gamma_{i_1}\cap\gamma_{i_2}\cap\gamma_{i_3}=\emptyset$ for $i_1,i_2,i_3$ all distinct.
\end{enumerate}
\end{corollary}
\begin{proof}
In view of Corollary \ref{networktriple} all of the junctions with multiplicity less than $2$ are triple junctions.
Let us assume that (i) is false and we will show that $V$ must satisfy (ii), that is, $V$ has no triple junctions so all of the singular points have multiplicity $2$.
If there is a geodesic segment of multiplicity $2$, then it cannot intersect any junction, because of the multiplicity bound.

The only possible junction is one formed by $4$ distinct geodesic segments of multiplicity one each.
We want to show that in this case it must be regular.
That is, at least two of the segments must have opposite directions at the singular point, which implies that so do the other two.

Denote by $v_1,v_2,v_3,v_4$ the unitary tangent direction of each geodesic segment at the singularity.
Let us suppose that at least $2$ of these are distinct and not opposite to each other. 
Without loss of generality we may assume it is $v_1$ and $v_2$.
Since we are in dimension $2$ we can use them as a basis and write $v_3$ and $v_4$ in terms of $v_1$ and $v_2$.
If one solves the system
\begin{equation*}
   \left\{
   \begin{aligned}
      & v_1+v_2+v_3+v_4 = 0 \,\text{(stationary condition)}\\
      & \|v_i\|=1 \,\text{(multiplicity one)}
   \end{aligned}
   \right.
\end{equation*}
then it is easy to see that, for example, $v_3$ must be opposite to either $v_1$ or $v_2$.

For the second part of $(ii)$, take $C\subset\supp\|V\|$ a connected component. 
If $V\restset C$ is given by a closed geodesic with multiplicity $2$ then the density condition implies that it cannot have junctions and the statement is true.
Otherwise, by what we showed above, each geodesic segment can be extended through the singular points.
Again, because of the density hypothesis we cannot have $3$ geodesics intersecting at the same point.
\end{proof}

The main result is a structure theorem for $1$-varifolds proved in \cite{wallard_falmgren1}. Here we state a particular case and refer to the original article for a proof.

\begin{theorem}\label{geodesicnetworkstructure}
Let $M$ be a closed manifold and $U\subset M$ an open set. If $V\in\intv_1(M)$ is stationary in $U$ then $V$ is a geodesic network in $U$.
\end{theorem}
\begin{proof}
Simply note that the definition of interval in \cite[\textsection 1]{wallard_falmgren1} is equivalent to being the image of a geodesic segment. The hypothesis for the theorem in \cite[\textsection 3]{wallard_falmgren1} are true because $V$ is integral. Finally, note that the set $S_V$ is the same as our set of junctions $\Sigma_V$.
\end{proof}

Now we prove the property that we are mainly interested for geodesic networks in $S^n$

\begin{prop}\label{massgeodesicnetwork}
Let $(S^n,g_0)$ be the round sphere of radius $1$, $U\subset S^n$ and $V\in\intv_1(S^n,U)$ be stationary in $U$ with total mass $\|V\|(S^n)<2\pi d$ for some positive integer $d$. Then $V$ is a geodesic network satisfying $\Theta^1(V,x)<d$ for all $x\in S^n$.
\end{prop}

\begin{proof}
We know by Theorem \ref{geodesicnetworkstructure} that $V$ is a geodesic network. Let us prove that $\|V\|(S^n)<2\pi d$ implies $\Theta^1(V,x)<d$ for every $x\in\supp\|V\|$.

Using proposition \ref{conemain}(\ref{conemaindensity}) and (\ref{conemainmass}), the Monotonicity formula for stationary varifolds (see \cite[\textsection 17.8]{lsimon1}) and $\unitvol_2=\pi$ we compute
\begin{equation*}
   \begin{aligned}
      \Theta^1(V,x) & = \Theta^2(C(V),x) = \lim_{r\rightarrow 0}\frac{\|C\|B(x,r)}{\unitvol_2 r^2} \\
                    & \leq \lim_{r\rightarrow\infty}\frac{\|C\|B(x,r)}{\unitvol_2 r^2} \\
                    & = \frac{\|V\|(S^n)}{2\unitvol_2} \\
                    & < d.
   \end{aligned}
\end{equation*}
\end{proof}

We can also prove a weaker version of this theorem for metrics that are sufficiently close to the round metric in $S^n$.
\begin{theorem}\label{densityellipsoid}
Let $g$ be a Riemannian metric in $S^n$. 
If $g$ is sufficiently $C^{\infty}$-close to the round metric then any varifold $W\in\intv_1(S^n)$ stationary with respect to the metric $g$ satisfying $\|W\|(S^n)<2\pi(d+\frac{1}{3})$ is a geodesic network such that $\Theta^1(W,x)\leq d$ for all $x\in\supp\|W\|$.
\end{theorem}
\begin{proof}
By Theorem \ref{geodesicnetworkstructure} we know that $W$ is a geodesic network with respect to the metric $g$.
It remains to prove the second statement.

Assume false, that is, there exist a sequence of metrics $g_i$ converging to $g_0$ and a sequence of integral varifolds $W_i$ stationary with respect to $g_i$ satisfying $\|W_i\|(S^n) < 2\pi(d+\frac{1}{3})$ and $\Theta^1(W_i,p_i)>d$ for some $p_i\in\supp\|W_i\|$.
In fact we must have $\Theta^1(W_i,p_i)\geq (d+\frac{1}{2})$ because $W_i$ is a geodesic network.

Since the first variation is continuous with respect to the metric, we may assume that each $W_i$ has bounded first variation in the metric $g_0$.
By the Compactness Theorem we may suppose that $W_i$ converges to an integral varifold $V$ stationary in the round metric and $p_i$ converges to $p\in\supp\|V\|$.

Furthermore, we have $\|V\|(S^n)\leq\liminf_{i\rightarrow\infty}\|W_i\|(S^n)<2\pi(d+\frac{1}{2})$ which implies that $\Theta^1(V,x)<(d+\frac{1}{2})$ for all $x\in\supp\|V\|$, by following a computation similar to the previous theorem. 
As before we must have $\Theta^1(V,x)\leq d$. 
On the other hand, the density is upper semicontinuous with respect to weak convergence of varifolds. 
In particular, $(d+\frac{1}{2}) \leq \Theta^1(W_i,p_i)\leq\Theta^1(V,p)\leq d$ for all $i$, which is a contradiction.

\end{proof}

%% file: wb-almostm-article2.tex
\section{Almost minimising varifolds}
In this section we define $\integer_2$-almost minimising varifolds and show that such $1$-dimensional varifolds cannot admit triple junctions.
%
%
\begin{definition}
Let $U\subset M$ be an open set, $\varepsilon>0$ and $\delta>0$. We define
\begin{equation*}
\acurr_k(U; \varepsilon,\delta)\subset\cycles_k(M)
\end{equation*}
as the set $T\in\cycles_k(M)$ such that any finite sequence $\{T_1,\ldots,T_m\}\subset\cycles_k(M)$ satisfying
\begin{enumerate}[(a)]
   \item $\supp(T-T_i)\subset U$ for all $i=1,2,\ldots m$;
   \item $\fnorm(T_i,T_{i-1})\leq\delta$ for all $i=1,2,\ldots m$ and
   \item $\mass(T_i)\leq\mass(T)+\delta$
\end{enumerate}
must also satisfy
\begin{equation*}
\mass(T_m)\geq\mass(T)-\varepsilon.
\end{equation*}
\end{definition}

Roughly speaking, if  $T$ belongs to $\acurr_k(U;\varepsilon,\delta)$ then any deformation of $T$ supported in $U$ that does not increase mass at least $\varepsilon$ must be $\delta$-far from $T$ in the $\fnorm$ metric.
Note that we define the elements of $\acurr_k$ as closed cycles in $M$ instead of relative cycles as defined in \cite{jpitts1}.
%
%
\begin{definition}
We say that a varifold $V\in\varf_k(M)$ is $\integer_2$-\textbf{almost minimising} in $U$ if for every $\varepsilon>0$ there exists $\delta>0$ and $T\in\acurr_k(U;\varepsilon,\delta)$ such that
\begin{equation*}
\fmetric(V,|T|)<\varepsilon.
\end{equation*}
\indent A varifold $V\in\varf_k(M)$ is said to be $\integer_2$-\textbf{almost minimising in annuli} if for every $p\in\supp \|V\|$ there exists $r>0$ such that $V$ is $\integer_2$-almost minimising in the annulus $A=A(p;s,r)$ for all $0<s<r$.
\end{definition}

\begin{definition}\label{no mass concentration}
For a cubical subcomplex $X\subset I^N$, we say that a flat continuous map $f:X\rightarrow \cycles_k(M)$ has \textit{no concentration of mass} if
\begin{equation*}
   \lim_{r\rightarrow 0}\sup\{\|f(x)\|(B(q,r)):x\in X \text{ and } q\in M\}=0.
\end{equation*}
\end{definition}

The next theorem shows the existence of $\integer_2$-almost minimising varifolds.
\begin{theorem}\label{existence almost minimising}
Let $X\subset I^N$ be a cubical subcomplex and $f:X\rightarrow \cycles_k(M)$ be a $p$-sweepout with no concentration of mass.
Denote $\Pi_f$ the class of all flat continuous maps $g:X\rightarrow \cycles_k(M)$ with no concentration of mass that are flat homotopic to $f$ and write
\begin{equation*}
   L[\Pi_f]=\inf_{g\in\Pi_f}\sup_{x\in X}\mass(g(x)).
\end{equation*}
If $L[\Pi_f]>0$ then there exists $V\in\intv_k(M)$ such that
\begin{enumerate}[(i)]
   \item $\|V\|(M)=L[\Pi_f]$;
   \item $V$ is stationary in $M$;
   \item $V$ is $\integer_2$-almost minimising in annuli.
\end{enumerate}
\end{theorem}

This was first proven by Pitts (see \cite[\textsection 4.10]{jpitts1}), for another proof (when $k=dim(M)-1$) we refer to \cite{fmarques-aneves1}.

Note that this is a weaker statement than in \cite{fmarques-aneves1}, but it remains true for all dimensions and codimensions.
This is because for every flat continuous homotopy class we can construct a discrete homotopy class just as in \cite[Theorem 3.9]{fmarques-aneves1} with the same width.
The final statement then follows from \cite[\textsection 4.10]{jpitts1}.


%
\begin{definition}
Let $T\in\cycles_k(M)$ and $W\subset M$ be an open set. We say that $T$ is locally mass minimising in $W$ if for every $p\in\supp(T)\bigcap W$ there exists $r_p>0$ such that $B(p,r_p)\subset W$ and for all $S\in\cycles_k(M)$ satisfying $\supp(T-S)\subset B(p,r_p)$ we have
\begin{equation*}
\mass(S)\geq\mass(T).
\end{equation*}
\end{definition}

\indent In the one dimensional case we have the following characterization:

%
\begin{prop}
Let $W\subset M$ be an open set, $Z\subset W$ compact and $T\in\cycles_1(M)$ be locally mass minimising in $W$. Then each connected component of $\supp(T)\bigcap Z$ is the restriction of a geodesic segment with endpoints in $W\setminus Z$.
\end{prop}
\begin{proof}
Let $A\subset\supp(T)\bigcap Z$ be a connected component. 
Cover $A$ by finitely many balls $B_i=B(p_i,r)$, $i=1,\ldots m$ such that each ball is contained in a convex neighborhood and $r<r_{p_i}$ for all $i$. 
Denote $C=\supp(T)\bigcap (B_1\cup\ldots\cup B_m)$, then each component $C\bigcap B_i$ is the unique minimising geodesic connecting the two points in $C\bigcap\partial B_i$. 
In particular the endpoints $A\bigcap \partial{Z}$ belong to the interior of a geodesic segment with endpoints in $\interior{Z}$ and $W\setminus Z$. We conclude that $A$ is given by the image of a broken geodesic with singular points in the interior of $Z$.\\
\indent Now, for each singular point $q\in A$ there exist $r_q$ such that $T\restset B(q,r_q)$ is mass minimising relative to its boundary. 
Thus it must be a geodesic segment, that is, $q$ is a smooth point in $A$. 
This implies that $C$ is the image of a geodesic segment with endpoints in $W\setminus Z$. 
The proof finishes by simply noting that $A=\closure{C\bigcap Z}$.
\end{proof}

\begin{corollary}\label{cyclestructure}
Let $W\subset M$ be an open set, $Z\subset W$ a compact set and $T\in\cycles_1(M)$ be locally mass minimising in $W$. Then, viewing $T$ as an integer coefficient current,
\begin{equation*}
T\restset Z = \sum_{i=1}^{k}\rectc(\beta_i,[1],\dot{\beta_i}),
\end{equation*}
where $\beta_i:[0,1]\rightarrow Z$ are geodesic segments for each $i=1,\ldots, k$ with endpoints in $\partial{Z}$.

In particular, the associated varifold $|T|\in\intv_1(M)$ is stationary in $W$.
\end{corollary}
\begin{proof}
We simply need to apply the Constancy Theorem (see \cite[\textsection 41]{lsimon1}) to each connected component. 
Since we are working with $\integer_2$ coefficients the density in each segment must be constant $1$.
\end{proof}

%
\indent The replacement theorem for almost minimising varifolds can be stated as follows:
\begin{theorem}\label{replacement}
Let $U\subset M$ be an open set, $K\subset U$ compact and $V\in\varf_k(M)$ $\integer_2$-almost minimising in $U$. There exists a non-empty set $\replacement(V;U,K)\subset\varf_k(M)$ such that every $V^*\in\replacement(V;U,K)$ satisfy:
\begin{enumerate}[(i)]
   \item $V^*\restset\grassm{k}{M\setminus K}=V\restset\grassm{k}{M\setminus K}$;
   \item $\|V^*\|(M)=\|V\|(M)$;
   \item $V^*$ is $\integer_2$-almost minimising in $U$;
   \item $V^*\restset\grassm{k}{\interior{K}}\in\intv_k(M)$ and
   \item for each $\varepsilon>0$ there exists $T\in\cycles_k(M)$ locally mass minimising in $\interior{K}$ such that $\fmetric(V^*,|T|)<\varepsilon$.
\end{enumerate}
\end{theorem}

\begin{proof}
The proof of (i)-(iv) is exactly as in \cite[\textsection 3.11]{jpitts1}. To show (v) one need to modify the construction in \cite[\textsection 3.10]{jpitts1} using our definition of almost minimising. 
\end{proof}

\begin{remark}\label{stationary replacement}
Note that if $V$ is stationary on all of $M$ then so is $V^*$.

In fact, $V^*$ is almost-minimising in $U$ (property \ref{replacement}(iii)) so it is also stationary in $U$.
Since $V^*$ coincides with $V$ on $M\setminus K$ then it is also stationary in $M\setminus K$.
That is, $V^*$ is stationary in $U$, $M\setminus K$ and $U\cap(M\setminus K)$.
Hence $V^*$ is stationary in $M$.
\end{remark}

\subsection{Almost minimising Geodesics Networks}
Here we will treat the particular case when $V$ is a geodesic network. 
Our main goal is to prove that the almost minimising property excludes the existence of triple junctions.

The rough idea is to use the replacement theorem and approximate $V^*$ by closed currents with coefficients in $\integer_2$. 
We will show that $V^*$ can be described as a non-zero $\integer_2$-cycle but triple junctions always have boundary in $\integer_2$.
From now on, given a varifold $V$ we will denote by $V^*$ a replacement given by Theorem \ref{replacement} whenever $V$ satisfies the conditions of the theorem.

\indent To prove the next technical lemma we will need the following theorem proven in \cite{bwhite3} by B.White and is used to prove a maximum principle for varifolds.

\begin{theorem}\label{maximumprinciplevector}
Let $N$ be a $n$-dimensional Riemannian manifold with boundary and $p\in\partial N$ such that $\kappa_1(p)+\ldots+\kappa_m(p) > \eta$, where $\kappa_1\leq\ldots\leq\kappa_{n-1}$ are the principal curvatures of $\partial N$ with respect to the inward normal vectorfield $\nu_N$. Then, given $\varepsilon>0$ there exists a supported vectorfield $X$ on $N$ such that $X(p)\neq 0$ is normal to $\partial N$ and
\begin{equation*}
   \langle X,\nu_N\rangle \geq 0 \text{ in } \partial N
\end{equation*}
and
\begin{equation*}
   \delta V(X)\leq -\eta\int|X|d\|V\|
\end{equation*}
for every $V\in\varf_m(N)$.
\end{theorem}

We remark that the same theorem is true with all its inequalities reversed, the proof is exactly the same (see \cite{bwhite3}).

\begin{corollary}\label{maximumprinciple1}
Let $M$ be a closed Riemannian manifold and $N$ an open set with strictly convex boundary with respect to the inward normal vectorfield ($\kappa_1\geq\eta>0$). 

If $V\in\varf_1(M)$ is stationary, $p\in\supp\|V\|\cap\partial N$ and $\supp\|V\|\cap B(p,\varepsilon)\cap N\neq\emptyset$ then $\supp\|V\|\cap B(p,\varepsilon)\cap (M\setminus\closure N) \neq\emptyset$.
\end{corollary}

\begin{proof}
First suppose there exists $\varepsilon>0$ such that $\supp\|V\|\cap B(p,\varepsilon)\subset \closure{N}$, that is, $W=V\restset\grassm{1}{B(p,\varepsilon)}\in\varf_1(\closure{N})$. 
Since $\partial N$ is strictly convex, we can choose $\eta>0$ in Theorem \ref{maximumprinciplevector} and obtain a vectorfield $X$ in $\closure{N}$ such that $\supp (X)\subset \closure N\cap B(p,\varepsilon)$ and
\begin{equation*}
   \delta W(X)+\frac{\eta}{2}\int|X|<0.
\end{equation*}
This is not a contradiction yet because $X$ is not a smooth vectorfield in $M$. However, we can construct a extension $\tilde X$ such that $\supp(\tilde{X})\subset B(p,\varepsilon)$, $\tilde X$ is $C^1$-close to $X$ and
\begin{equation*}
   \delta W(\tilde X)+\frac{\eta}{2}\int|\tilde X|<0.
\end{equation*}
By construction $\supp(\tilde X)\subset B(p,\varepsilon)$ hence $\delta V(\tilde X)=\delta W(\tilde X)<0$. This is a contradiction because $V$ is stationary, thus $\supp\|V\|\cap B(p,\varepsilon)\cap (M\setminus\closure N) \neq\emptyset$.
\end{proof}

We now show that an almost-minimising geodesic network is its own replacement.
To simplify notation, from now on we write $V\restset U= V\restset\grassm{1}{U}$ whenever $V\in\intv_1(M)$ and $U\subset M$ is an open set.

\begin{lemma}\label{uniquereplacement}
Let $M$ be a closed Riemannian manifold and $V\in \intv_1(M)$ be a geodesic network and $p\in\Sigma_V$ be a junction point. 
If $V$ is almost minimising in annuli at $p$ then there exists $r>0$ and a compact set $K\subset A(p;r,3r)$ such that
\begin{enumerate}[(i)]
   \item $V$ is almost minimising in $A(p;r,3r)$ and
   \item $\replacement(V; A,K)=\{V\}$.
\end{enumerate}
\end{lemma}
\begin{proof}
Since $V$ is a geodesic network, then its singularities are isolated.
That is, there exists $r_p>0$ such that $p$ is the only singularity in $B(p,r_p)$.

Firstly choose $r>0$ such that $4r < r_p$, $B = B(p,4r)$ is a convex ball and $V$ is almost minimising in $A=A(p;r,3r)$.
It follows from the structure of a geodesic network that
\begin{equation*}
V\restset{B}=\sum_{j=1}^m\rectv(\alpha_j,\theta_j)
\end{equation*}
where $\alpha_j:[0,4r]\rightarrow\closure B$ is a minimising geodesic parametrized by arc-length for each $j=1,\ldots,m$. 
By abuse of notation we identify the curves $\alpha_j$ with its image.

Secondly, we can choose $\delta<r$ sufficiently small such that the balls $K_j=\bar B(\alpha_j(2r),\delta)$ have strictly convex boundary with respect to the inward normal vector and are pairwise disjoint.
Define $a_j=\alpha_j(2r-\delta)$, $b_j=\alpha_j(2r+\delta)$ and $K=K_1\cup\ldots\cup K_m\subset A$.

Finally we take $V^*\in\replacement(V;A,K)$ a replacement for $V$ and define $V^*_j=V^*\restset{\interior{K_j}}$ and $V_j=V\restset{\interior{K_j}}$.
By property \ref{replacement}(i) it is sufficient to show that $V^*_j=V_j$ for each $j$.

\begin{claim}
$\sum_{j=1}^m \|V^*_j\|(M)=\sum_{j=1}^m\|V_j\|(M)$
\end{claim}
This follows directly from properties \ref{replacement}(i) and (ii).

\begin{claim}
For each $j=1,\ldots,m$ either $V^*_j=0$ or $\supp\|V^*_j\|$ contains a rectifiable curve connecting $a_j$ to $b_j$.
\end{claim}
Note that $\supp\|V^*_j\|$ only intersects $\bdry{K_j}$ at the points $a_j$ and $b_j$.
In fact, suppose there is another point of intersection.
Then, by the maximum principle (Corollary \ref{maximumprinciple1}) it follows that $\supp\|V^*\|\setminus\interior{K_j}=\supp\|V^*\restset{M\setminus K}\|$ also contains that point, but this contradicts property \ref{replacement}(i).

Now, suppose $\supp\|V^*_j\|$ contains no curve joining $a_j$ and $b_j$.
In that case, we can write $\supp\|V^*_j\|=C_a\cup C_b$ where $C_a$ and $C_b$ are closed disjoint sets containing $a_j$ and $b_j$ respectively (these are not unique and not necessarily connected).
Take $U_a$ and $U_b$ open and disjoint neighbourhoods of $C_a$ and $C_b$ in the interior of $K_j$ respectively.
We will show that $V_j^*\restset{U_a}=V_j^*\restset{U_b}=0$.

Take for example $V_j^*\restset{U_a}$, which is stationary (see remark after Theorem \ref{replacement}).
Now, consider $B(\sigma)=B(\alpha(2r\sigma),\sigma\delta)$ then $V_j^*\restset{U_a}$ is entirely contained in $B(1)=\interior{K_j}$ and it only intersects the boundary at the point $a_j$.
Since $\partial B(\sigma)$ is strictly convex for all $\sigma$, a maximum principle argument shows that $V_j^*\restset{U_a}$ is contained in $B(\sigma)$ for all $\sigma<1$ thus proving that $V_j^*\restset{U_a}=0$.
The same argument shows that $V_j^*\restset{U_b}=0$ and we prove the claim.

\input{wb-convex}
%

\begin{claim}
$V_j^*\neq 0$ for all $j=1,\ldots,m$
\end{claim}
Consider $B_j'=B(\alpha_j(2r),\delta')$ with $\delta<\delta'<r$ such that $K_j\subset B'_j\subset A$ are still pairwise disjoint.
Then property \ref{replacement}(i) implies that
\begin{equation*}
   V^*\restset{B_j'\setminus K_j}=\rectv(\alpha_j\cap(B_j'\setminus K_j),\theta_j).
\end{equation*}
If $V^*_j$ was zero, then in particular $V^*_j=\rectv(\alpha_j,0)$.
But $V^*\restset{B_j'}=V^*\restset{B_j'\setminus K_j}+V_j^*$ is stationary and its support is contained in $\alpha_j$.
From the Constancy Theorem we conclude that $\theta_j=0$ which is a contradiction, thus $V_j^*\neq 0$.

This means that $\supp\|V_j^*\|$ contains a rectifiable curve $C_j$ connecting $a_j$ to $b_j$ for all $j=1,\ldots,m$.
In particular this implies that $l(C_j)\geq d(a_j,b_j)$.
Since $V_j^*$ is integral (see property \ref{replacement}(iv)) it follows that $\|V_j^*\|(M)\geq d(a_j,b_j)$.
However, $\alpha_j\cap K_j$ is a minimising geodesic connecting $a_j$ to $b_j$, so $\|V_j\|(M)=d(a_j,b_j)$.
We conclude that $\|V^*_j\|(M)\geq \|V_j\|(M)$ for all $j=1,\ldots,m$. 
Claim $1$ implies that we have in fact 
\begin{equation*}
   \|V^*_j\|(M)= \|V_j\|(M) \,\text{ for all }\, j=1,\ldots,m.
\end{equation*}

On the other hand, we have $d(a_j,b_j)=\|V_j^*\|(M)\geq l(C_j) \geq d(a_j,b_j)$, that is, $C_j$ is a minimising curve and it must be a geodesic.
Since $\alpha_j\cap K_j$ is the unique geodesic connecting $a_j$ to $b_j$ we conclude that $C_j=\alpha_j\cap K_j$.
Finally, this implies that $\supp\|V_j^*\|=\supp\|V_j\|$ because otherwise there would be more contribution of mass.
Applying the Constancy Theorem again we show that $V_j^*=V_j$ and this finishes the proof.

\end{proof}

The last result we need relates flat convergence of $\integer_2$-currents and the weak convergence of the associated varifold. This was proven in \cite{bwhite2} by B.White.
\begin{theorem}\label{flatversusweakconvergence}
Let $M$ be a Riemannian manifold, $\{W_i\}\subset\intv_k(M)$ be a sequence converging to an integral varifold $W$. Suppose that:
\begin{enumerate}[(I)]
   \item Each $W_i$ has locally bounded first variation;
   \item $\partial[W_i]$ converges in the flat topology.
\end{enumerate}
Then $[W_i]$ converges to $[W]$ in the flat topology.
\end{theorem}

Finally we prove our main result of this section.

\begin{theorem}\label{integerdensity}
Let $M$ be a closed surface, $V\in\intv_1(M)$ a geodesic network and $p\in\Sigma_V$ a junction point. 
If $V$ is $\integer_2$-almost minimising in annuli at $p$, then 
\begin{equation*}
\Theta^1(V,p)\in\naturals.
\end{equation*}
In particular $p$ is not a triple junction.
\end{theorem}

\begin{proof}
Let $r>0$, $B=B(p,4r)$, $A=A(p;r,3r)$ and $K\subset A$ as in Lemma \ref{uniquereplacement}. 
Applying property \ref{replacement}(v), Corollary \ref{cyclestructure} and the Compactness theorem for $\integer_2$-chains (see \cite[Theorem 5.1]{bwhite2}) we may assume there exists a convergent sequence $\{T_i\}_{i\in\naturals}\subset\cycles_1(M)$ and $T\in\cycles_1(M)$ such that
\begin{enumerate}[(a)]
   \item $T_i\rightarrow T$ in the $\fnorm$-norm;
   \item $V_i=|T_i|$ is stationary in $\interior{K}$ and
   \item $V_i\rightarrow V$ in the $\fmetric$-metric.
\end{enumerate}

Even though convergence of chains in the flat norm do not correspond to weak convergence for varifolds, in the stationary case, with convergent boundary, it does.

We want to apply Theorem \ref{flatversusweakconvergence} for the sequence $\{V_i\restset{\interior{K}}\}_{i\geq 1}$. 
We know that $\partial[V_i\restset{\interior{K}}]\rightarrow\partial T\restset\interior{K}$ by the definition of $V_i$.
Together with property (b) it means that the sequence satisfies the hypothesis of the theorem.
We conclude that
\begin{equation*}
   T\restset\interior{K} = [V\restset{\interior{K}}].
\end{equation*}

Since $V\restset{B}=\sum_{j=1}^m\rectv(\alpha_j,\theta_j)$ for some geodesic segments $\alpha_j$ and $\theta_j\in\integer_{>0}$, we have
\begin{equation*}
[V\restset{\interior{K}}]=\sum_{j=1}^m\rectv(\alpha_j\cap{\interior{K}},[\theta_j])
\end{equation*}
and $[\theta_j]$ is non-zero only when $\theta_j$ is odd. 

If $\theta_j$ is even for all $j$ then the density at $p$ must be an integer and we finish our proof because geodesic segments with even multiplicity contribute to the density at $p$ with an integer number.\\
\indent In case some $\theta_j$ is odd we have that $T\neq 0$ and $\supp(T)\subset\supp\|V\|$.
We can view $T\restset B$ as an integer chain and apply the Constancy theorem for integral currents (see \cite[\textsection 26.27]{lsimon1}) and the fact that $T$ and $[V]$ coincide in $\interior{K}$ to conclude that
\begin{equation*}
T\restset B=\sum_{j=1}^m\rectc(\alpha_j,[\theta_j],\dot{\alpha_j}).
\end{equation*}

Now we simply note that $p$ is a boundary point for $T$ unless the number of $\theta_j$ such that $[\theta_j]\neq 0$ is even.
That is, there is an even number of geodesic segments $\alpha_j$ with odd multiplicity and in particular its density contribution is an integer number.
This finishes the proof because $T$ is a closed chain.

\end{proof}

%% file: wb-convex.tex


\newcommand{\centrearc}[4]{\draw #1 ++ (#2:#4) arc (#2:#3:#4)}


\begin{center}
\begin{tikzpicture} 
[point/.style = {circle,fill,inner sep=0pt,minimum size=3pt,draw},
labeling/.style = {inner sep=1pt,minimum size=1pt},
labeldist/.style = {circle},
scale = 1.3,
every node/.style={scale=0.8}
]

   \node[point,label={[label distance = -0.2 cm]150:$p$}] (p) at (4,0) {};
   \draw (p) -- (-4.1,0);
   \draw (p)++(60:0.3 cm) -- (p);
   \draw (p)++(-60:0.3 cm) -- (p);
   \node[circle,label={[label distance = -0.3 cm]150:$\alpha_j$}] (alphaj) at (3,0) {} ;
   \foreach \i in {2,3,4}{
   \pgfmathtruncatemacro{\pos}{4-2*\i}
      \draw (\pos,-1pt)--(\pos,1pt) node[labeling,label=-150:$\i r$] (\i r) {};
};
   \draw (2,-1pt)--(2,1pt) node[below right] (r) {$r$};
   \draw[line width = 1pt] (p)++(185:8 cm) arc (185:175:8cm); 
   \node[labeling,label={[label distance=-0.1 cm]120:$B$}] at ($(175:8 cm)+(p)$) {};
   \draw[thin, dashed] (p)++(195:6 cm) arc (195:165:6cm); 
   \node[labeling,label={[label distance=-0.1 cm]120:$A$}] at ($(165:6 cm)+(p)$) {};
   \draw[thin, dashed] (p)++(195:2 cm) arc (195:165:2cm); 
   \draw[line width=1pt] (2r) circle (1.6 cm); 
   \node[label=90:$K_j$] at (30:1.8 cm) {};
   \node[label={[label distance = -0.2 cm]-45:$a_j$}] at (1.6,0) {};
   \node[label={[label distance = -0.3 cm]225:$b_j$}] at (-1.6,0) {};
   \draw[dotted] (0.5,0) circle (1.1 cm);
   \node[label=50:$B(\sigma)$] at (150:1.1 cm) {};
   \draw [densely dashed] (1.6,0) .. controls +(120:3 cm) and +(240:3 cm) .. (1.6,0);
   \draw [densely dashed] (2r)++(0.3,0.5) circle (0.2 cm);
   \node[label={[label distance = 0.1 cm]80:$C_a$}] at (0.4,0.5) {};
\end{tikzpicture}

\end{center}


%% file: wb-ellip-article.tex
\section{The width of an Ellipsoid}

Here we will apply the previous results to estimate some of the $k$-width of ellipsoids sufficiently close to the round sphere.

\subsection{Sweepouts of $S^2$}
Let $(S^2,g_0)$ denote the round $2$-dimensional sphere with radius $1$ in $\real^3$.
We will construct $k$-sweepouts of $S^2$ as families of algebraic sets in $\real^3$.
This is similar to how it is done in \cite{lguth1} for the unit ball.

Denote by $x_1,x_2,x_3$ the coordinates in $\real^3$ with respect to the standard basis.
Let $p_i:\real^3\rightarrow \real$ denote the following polynomials for $i=1,\ldots,8$:
\begin{equation*}
   \begin{aligned}
      p_j(x) & = x_j \text{ for } j=1,2,3;\\
      p_4(x) & = {x_1}^2;\\
      p_5(x) & = x_1 x_2;\\
      p_6(x) & = x_1 x_3;\\
      p_7(x) & = x_2 x_3;\\
      p_8(x) & = {x_3}^2.
   \end{aligned}
\end{equation*}
Note that we skipped the polynomial ${x_2}^2$. 
The reason for this is because we are only interested in the zero set restricted to the sphere, which is given by the equation ${x_1}^2+{x_2}^2+{x_3}^2-1=0$.
That is, $p_4$, $p_8$ and $x\mapsto {x_2}^2$ are linearly dependent so their linear combinations will define the same algebraic sets.

Now, put $A_k=\span_{\real}\left({1}\cup_{j=1}^{k}{p_j}\right)\setminus\{0\}$ and define the relation $q\sim\lambda q$, for $\lambda>0$ and $q\in A_k$.
Note that the zero set is invariant under this relation, that is, $\{\lambda q=0\}=\{q=0\}$ so it makes sense to define the map $F_k:(A_k/\sim)\rightarrow\cycles_1(S^2)$ as
\begin{equation*}
   F_k([q])=\partial[\{q \leq 0\}\cap S^2],
\end{equation*}
where $[R]$ denotes the mod $2$ current associated with $R\subset S^2$.
It is clear that $F_k$ is well defined and it takes values in $\cycles_1(S^2)$.

We can identify $(A_k/\sim)$ with $\real P^{k}$ and we observe that $F_k$ is flat continuous and it defines a $k$-sweepout.
The proof is exactly the same as in \cite[\textsection 6]{lguth1} with the appropriate adaptations and we leave it to the reader.

\begin{lemma}
Let $F_k:\real P^k\rightarrow \cycles_1(S^2)$, k=1,\ldots,8, be the family of cycles defined above.
Then $F_k$ has no concentration of mass.
\end{lemma}
\begin{proof}
Take $p\in S^2$ and $0<r<\pi$ and denote by $\alpha_p$ the equator given by $p^{\perp}\cap S^2$, where $p^{\perp}$ is the plane normal to $p$ in $\real^3$.
Consider the ball $B(p,r)\subset S^2$.
We can parametrize the space of geodesics that go through $B(p,r)$ as $G(r)=\{q^{\perp}\cap S^2 : d(q,\alpha_p)<r\}$. The set $G(r)$ defines a spherical segment whose area is $area(G(r))=4\pi\sin(r)$.

If $x\in\real P^k$ is such that $F_k(x)\cap B(p,r)\neq\emptyset$ then it follows from the Crofton formula that
\begin{equation*}
\mass(F_k(x)\restset B(p,r))=\frac{1}{4}\int_{\Gamma\in G(r)}\#(\Gamma\cap F_k(x)).
\end{equation*}
Since $\Gamma\cap F_k(x)$ is the intersection of a plane with $S^2$ and $F_k(x)$ then it is the solution of a system of $3$ polynomials of degree $1$, $2$ and at most $2$ ($1$ if $k=1,2,3$ or $2$ if $k=4,\ldots, 8$), respectively.
It follows that the intersection is generically $\#(\Gamma\cap F_k(x))\leq 4$.
Hence,
\begin{equation*}
\mass(F_k(x)\restset B(p,r))\leq 4\pi\sin(r).
\end{equation*}
If we take $r\rightarrow 0$ we conclude that $F_k$ has no concentration of mass at $p$.
Since $p$ was arbitrary we conclude the proof.
\end{proof}

\begin{remark}
Note that the same proof is valid for any family of algebraic curves in $S^2$ with bounded degree.
\end{remark}


\begin{theorem}\label{width sphere}
If $S^2$ is the round $2$-sphere of radius $1$, then
\begin{enumerate}[(i)]
   \item $\width_1(S^2)=\width_2(S^2)=\width_3(S^2)=2\pi$;
   \item $\width_4(S^2)=\width_5(S^2)=\width_6(S^2)=\width_7(S^2)=\width_8(S^2)=4\pi$.
\end{enumerate}
\end{theorem}
\begin{proof}
\textbf{(i):} By the Crofton formula we have that, $\mass(F_k(q))\leq 2\pi$ for all $q\in\real P^k$ and $k=1,2,3$.
In fact, it is not hard to see that $\sup \mass(F_k(q))=2\pi$.
That is, $\width_k\leq 2\pi$.

Suppose $\width_k< 2\pi$, then there exists another $k$-sweepout with no concentration of mass $\tilde{F}$ such that $L[\Pi_{\tilde{F}}]<2\pi$. 
Hence, Theorem \ref{existence almost minimising} would give us a stationary $\integer_2$-almost minimising integral varifold with $\|V\|(S^2)<2\pi$.
This is a contradiction because Theorem \ref{massgeodesicnetwork} tells us that the density would be lower than $1$ everywhere.
So $F_k$ is optimal and $\width_k=2\pi$ for $k=1,2,3$.

For the next item we need a lemma whose proof we give in the Appendix.

\begin{lemma}\label{increasing width}
Let $S^2$ be the round $2$-sphere of radius $1$, then $\width_4>2\pi$.
\end{lemma}

\noindent \textbf{(ii):} When $k=4,5,6,7,8$ the degree of the polynomials are less than or equal to $2$, thus, using the Crofton formula again, $\mass(F_k(q))\leq 4\pi$ for all $q\in\real P^k$.
As before, it is trivial to check that $\sup\mass(F_k(q))=4\pi$ from which we get $\width_k\leq 4\pi$.

By Lemma \ref{increasing width} and the previous item we already know that $\width_k\geq \width_4>2\pi$.
Suppose $\width_k<4\pi$ then, as before, we have a $k$-sweepout $\tilde{F}$ with no concentration of mass such that $\width_k\leq L[\Pi_{\tilde{F}}]<4\pi$.
From Theorem \ref{existence almost minimising} we produce $V\in\intv_1(S^2)$ stationary and $\integer_2$-almost minimising.
It follows from Theorems \ref{geodesicnetworkstructure} and \ref{integerdensity} that $V$ has density constant to $1$.
Hence $V$ corresponds to a closed regular geodesic, that is, $\|V\|(S^2)=2\pi$, which is a contradiction.
\end{proof}

\subsection{Geodesics on Ellipsoids}

Our goal here is to find the varifold that realizes the $k$-width of an ellipsoid sufficiently close to the round sphere.

Let $E^2=E^2(a_1,a_2,a_3)$ be an ellipsoid defined by the equation $a_1x_1^2+a_2x_2^2+a_3x_3^2-1=0$ in $\real^3$.
If the parameters $a_1,a_2,a_3$ are all sufficiently close to $1$ then it is clear that the induced metric in $E^2$ is $C^{\infty}$-close to the round metric in $S^2$.
We can assume other properties that we summarize here.
\begin{prop}\label{genericellipsoid}
Let $\gamma_i=\{x_i=0\}\cap E^2$ for $i=1,2,3$ be the three principal geodesics in $E^2$, $\gamma_i^{(r)}$ be the $r$-covering of $\gamma_i$ for $r\in\naturals$ and $\width_k(E^2)$ denote the $k$-width for $k\in\naturals$.
If we choose $a_1<a_2<a_3$ sufficiently close to 1 then the following is true:
\begin{enumerate}[(i)]
   \item $2\pi(1-\frac{1}{4})<L(\gamma_1)<L(\gamma_2)<L(\gamma_3)<2\pi(1+\frac{1}{4})$; 
   \item $\index(\gamma_i^{(r)})=i+2(r-1)$ and $\nullity(\gamma_i^{(r)})=0$ for $i=1,2,3$ and $r<100$;
   \item if $\alpha$ is a smooth closed geodesic with $L(\alpha)<100\pi$ then $\alpha = \gamma_i^{(r)}$ for some $i=1,2$ or $3$ and $r>0$;
   \item $|\width_k(E^2)-\width_k(S^2)|<\frac{1}{4}$ for all $k<100$;
\end{enumerate}
By $\index(\gamma)$ and $\nullity(\gamma)$ we mean the Morse index and nullity as smooth closed geodesics, that is, critical points of the energy functional.
\end{prop}
\begin{proof}
\textbf{(i)}: Note, for example, that $\gamma_1$ is a planar ellipsis with axes $\frac{1}{a_2}$ and $\frac{1}{a_3}$, similarly for the other two.
So, as long as $a_i$ are close to $1$ each ellipsis is close to a circle of length $2\pi$.

\indent \textbf{(ii)}: See \cite[XI,Theorem 3.3]{mmorse1}.

\indent \textbf{(iii)}: See \cite[XI,Theorem 4.1]{mmorse1}.

\indent \textbf{(iv)}: Note that every sweepout of $E^2$ is also a sweepout for $S^2$, simply by the fact they are both diffeomorphic and the definition of sweepout does not depend on the metric. Since the metric in $E^2$ can be chosen sufficently close to the round metric we can prove that each $k$-width is continuous by simply using the same approximating sweepouts. The uniform convergence follows directly because we are considering only finitely many $k$-widths.
\end{proof}

Given these three main ellipses we are able to define the varifolds that will be candidates to realize the first $8$ widths of $E^2$. 
Define
\begin{equation*}
   \begin{aligned}
      W_j & = \rectv(\gamma_j,1), j=1,2,3;\\
      W_4 & = \rectv(\gamma_1,2);\\
      W_5 & = \rectv(\gamma_1,1)+\rectv(\gamma_2,1);\\
      W_6 & = \rectv(\gamma_2,2);\\
      W_7 & = \rectv(\gamma_1,1)+\rectv(\gamma_3,1);\\
      W_8 & = \rectv(\gamma_2,1)+\rectv(\gamma_3,1);\\
      W_9 & = \rectv(\gamma_3,2).
   \end{aligned}
\end{equation*}
\begin{remark}
Suppose $E^2$ is sufficiently close to the round sphere of radius $1$.
Since these are all possible combinations of the three principal geodesic networks with density less than or equal to $2$, Theorems \ref{densityellipsoid}, \ref{integerdensity} and Corollary \ref{characterization} imply that these are also the only almost minimising geodesic networks with mass less than $2\pi(2+\frac{1}{3})$.

They also correspond to the zero set (counted with multiplicity) of the polynomials $p_j$, defined in the previous section, intersected with $E^2$ (except for $W_6$).
\end{remark}

Before proceeding to the main theorem we need a technical lemma that was proved in \cite[\textsection 6]{fmarques-aneves1} under a different context.
We explain how to obtain our result from their proof in the Appendix.
\begin{lemma}\label{ellipsoid width increasing}
Let $E^2$ be an ellipsoid as in Proposition \ref{genericellipsoid}. Then $\width_i<\width_{i+1}$ for $i=1,\ldots, 7$.
\end{lemma}
\begin{theorem}\label{ellipsoidwidth}
Let $E^2$ be an ellipsoid as in Proposition \ref{genericellipsoid}. The following holds:
\begin{enumerate}[(i)]
   \item if $i=1,2$ or $3$ then $\width_i(E^2)=\|W_i\|(E^2)$;
   \item if $i=4,\ldots 8$ then $\width_i(E^2)=\|W_l\|(E^2)$ for some $l=4,\ldots,9$ without repetition.
\end{enumerate}
\end{theorem}
\begin{proof}
Firstly, it follows from Proposition \ref{genericellipsoid}(iv) and Theorem \ref{width sphere}(i) that
\begin{equation*}
   \begin{aligned}
      \width_i(E^2)<2\pi(1+\frac{1}{4}) & \text{, for }i=1,2,3 \text{ and}\\
      \width_i(E^2)<2\pi(2+\frac{1}{4}) & \text{, for }i=4,\ldots,8.
   \end{aligned}
\end{equation*}

In either case we claim that there exists an optimal sweepout for $\width_i$.
Indeed, if no such map existed for some $i$ we would have a sequence of sweepouts $\{F_k\}$ satisfying $\width_i<L[F_{k+1}]<L[F_k]<2\pi(2+\frac{1}{4})$.
Each $F_k$ provides us a distinct almost-minimising geodesic network with mass less than $2\pi(2+\frac{1}{4})$ by Theorem \ref{existence almost minimising}, and the characterization of stationary integral varifolds (Theorem \ref{geodesicnetworkstructure}).
However, as we have already remarked, there only finitely many such varifolds (that is to say, the previously defined $W_j$) so we have a contradiction.

Secondly, Lemma \ref{ellipsoid width increasing} tells us that $\width_1<\ldots<\width_8$.
Hence, each optimal sweepout gives us an almost-minimising geodesic network $V_i$ satisfying $\|V_i\|(E^2)=\width_i(E^2)$ and $\|V_i\|(E^2)<\|V_{i+1}\|(E^2)$ for $i=1,\ldots,8$.

\noindent \textbf{(i):} 
For $i=1,2,3$ we have $\|V_i\|(E^2)<2\pi(1+\frac{1}{4})$ so each one of these must correspond to one $W_j$, $j=1,2,3$.
Since their masses are ordered as $\|W_1\|(E^2)<\|W_2\|(E^2)<\|W_3\|(E^2)$ we must have $V_i=W_i$, $i=1,2,3$.

\noindent \textbf{(ii):} 
For $j=4,\ldots,8$ the $W_j$'s are not necessarily ordered by their mass.
To be specific, we cannot guarantee for a general ellipsoid that $\|W_6\|(E^2)<\|W_7\|(E^2)$ or vice-versa.
However, we know that each $V_i$ corresponds to one of the $W_j$'s and this correspondence must be one to one, which finishes the proof.
\end{proof}

At last we give a counterexample to Question 1.

\begin{corollary}
Let $E^2$ be as in Proposition \ref{genericellipsoid}, then \textbf{\underline{Question 1}} is false for $E^2$.
\end{corollary}

\begin{proof}
First of all we observe that if the support of $V$ is given by a smooth closed geodesic $\gamma$ then $\index(V)$ and $\nullity(V)$ as a varifold are the same as $\index(\gamma)$ and $\nullity(\gamma)$ has a critical point for the energy functional.

Now, Theorem \ref{ellipsoidwidth}(ii) tells us that $\width_4(E^2)=\|W_j\|(E^2)$ for some $j=4,\ldots,9$. 
Where $W_j$, $j=1,\ldots,8$ are as before.
Since there are $6$ varifolds to choose for $5$ widths we know that one, and only one, will not correspond to a width.

The first $3$ varifolds are ordered as $\|W_4\|(E^2)<\|W_5\|(E^2)<\|W_6\|(E^2)$ which implies that $\width_4$ must correspond to either $W_4$ or $W_5$.
If $\width_4(E^2)=\|W_4\|(E^2)$ we are done because the number of parameters is $4$ and $\index(V)+\nullity(V)=3+0<4$, as given by property \ref{genericellipsoid}(ii).

If this is not the case, then $W_4$ is the only varifold that does not correspond to any width and all the other ones must correspond to one, and only one, width.
As we have already pointed out, the comparison between $\|W_6\|(E^2)$ and $\|W_7\|(E^2)$ is not known in general.
In any case, $W_6$ must correspond to either $\width_5$ (if $\|W_6\|(E^2)<\|W_7\|(E^2)$)  or $\width_6$ (if $\|W_7\|(E^2)<\|W_6\|(E^2)$).
On the other hand, $\index(W_6)+\nullity(W_6)=4+0 < 5 $, which disproves the conjecture in either case.
\end{proof}

\begin{remark}
The proof above also gives us an example of a unstable min-max critical $1$-varifold with multiplicity and smooth embedded support.
As described above, we will either have $\gamma_1$ or $\gamma_2$ (that have index $1$ and $2$ respectively) with multiplicity $2$.

It is conjectured (see \cite[p.2]{fmarques-aneves3}) that this should not happen for min-max critical hypersurfaces.
The main difference is because in the hypersurface case one could be able to de-singularize two minimal surfaces (for example two great spheres in $S^3$ approaching a sphere with multiplicity $2$) along their intersection and obtain a min-max ``competitor'' with very similar area, but embedded and with different topology.
Such procedure doesn't exist for curves and our example settles this question for the one dimensional case.
\end{remark}

%% file: wb-problems-article.tex
\section{Further Problems}

We would like to propose a general formula for the width of the round sphere $S^2$.
First let us give our conjecture and then explain the motivation.
We expect that
\begin{equation*}
   \width_j=2\pi k, \text{ if } j\in\{k^2,\ldots,(k+1)^2-1\}.
\end{equation*}
A simple computation shows that this would imply the Weyl law for $S^2$.
Of course, to prove the Weyl law it is not necessary to compute the width spectrum, one is only interested in its asymptotic behavior.
This is a much stronger conjecture.

Denote by $P[\real^3, d]$ the space of real polynomials of degree less than or equal to $d$ in $3$ variables.
For each $p\in P$ we can define $\{p=0\}\cap S^2$ as we have already done.
However, any polynomial that contains the fact $({x_1}^2+{x_2}^2+{x_3}^2-1)$ do not define a $1$-cycle in $S^2$ so we have to quotient these out.
That is, we are interested in the space $A_d=P[\real^3,d]/\langle{x_1}^2+{x_2}^2+{x_3}^2-1\rangle_d$, where $\langle{x_1}^2+{x_2}^2+{x_3}^2-1\rangle_d$ denotes the ideal generated by $({x_1}^2+{x_2}^2+{x_3}^2-1)$ intersected with $P[\real^3,d]$.

Now, note that we can write $A_{k}=A_{k-1}\oplus H_{k}$ where $H_{k}$ is the space of homogeneous polynomials in $3$ variables of degree $k$.
The space $H_{k}$ is isomorphic to the eigenspace of the ${k}^{th}$ eigenvalue of the Laplacian in $S^2$ and its dimension is $2k+1=(k+1)^2-k^2$.
Using the polynomials in $A_{k-1}$ and a basis for $H_k$ we can construct $j$-sweepouts for $j=k^2,\ldots,(k+1)^2-1$ whose minmax values are $2\pi k$ as given by the Crofton formula.
We expect these sweepouts to be optimal for the round sphere.

This is motivated by Lusternik-Schnirelmann theory on manifolds.
We believe that the width will be realised by a combination of great circles with possible multiplicities.
Lusternik-Schnirelmann theory indicates that if $\width_k=\width_{k+N}$ then there exists a $N$-parameter family of varifolds with constant mass $\width_k$.
More generally, one would expect the space of critical varifolds with mass $\width_k$ to have Lusternik-Schnirelmann category greater or equal to $N$
In the case of $S^2$ the space of $k$-combinations of great circles is simply the space of unordered $k$-tuples of great circles, that is, it is given by $SP^k(\real P^2)$.
We denote by $SP^k(X)$ the quotient of $X^k$ by the action of the $k$-symmetry group $S_k$.
It is known that $SP^k(\real P^2)=\real P^{2k}$ (see \cite{varnold1}), whose Lusternik-Schnirelmann category is $2k+1$.
Finally our conjecture implies that the equality gaps in the width spectrum are given by $\width_{k^2}=\width_{(k+1)^2-1}$, which is consistent with the Lusternik-Schnirelmann motivation.
As a brief remark we would like to point out the for higher dimensions the same ideas would violate the category of the critical set.

Unfortunately none of this has been proved.
Neither the category ideas or the optimality of the polynomial sweepouts are known.
The Lusternik-Schnirelmann theory for smooth functions on manifolds (see \cite{ocorneaetal1}) does not carry over to our case directly.

%% file: wb-appendix-article.tex
\begin{appendices}

\section{}

First let us extract a weaker version of the results in \cite{fmarques-aneves1}.
From the proof of \cite[Theorem 6.1]{fmarques-aneves1} we can obtain the following general, but weaker, proposition.

\begin{prop}\label{equality width}
Let $M$ be a closed Riemannian manifold and $\{\width_k(M)\}_{k\in\naturals}$ be the width spectrum corresponding to $1$-cycles in $M$. 
If $\width_k(M)=\width_{k+1}(M)$ for some $k$, then there exist infinitely many geodesic networks with mass $\width_k(M)$ and that are almost-minimising in annuli at every point.
\end{prop}

The proof is similar to \cite[Theorem 6.1]{fmarques-aneves1}, however we cannot use Schoen-Simon's Regularity Theorem or the Constancy Theorem (as in \cite[Claim 6.2]{fmarques-aneves1}).
To overcome this one notes that if a sequence of varifolds converge to a geodesic network then the sequence of associated currents converge to a subnetwork of the limit.

More precisely, let $\{T_i\}_{i\in\naturals}\subset \cycles_1(M)$ be a sequence of flat cycles such that $|T_i|\rightarrow V$ and $T_i\rightarrow T$.
If $V$ is a geodesic network defined by geodesic segments $\{\gamma_1,\ldots,\gamma_m\}$ and its respective multiplicities, then $T$ is a cycle (not necessarily stationary) defined by a subset of geodesics $\Omega\subset \{\gamma_1,\ldots,\gamma_m\}$ with multiplicity one each.

This is true because the support of the limit is contained in the varifold geodesic network, then we can apply the Constancy Theorem to each geodesic segment whose intersection is non-empty.
If we assume that the set of geodesic networks is finite, then so is the set of all possible subnetworks (not necessarily stationary) and the rest of the proof is the same as in \cite{fmarques-aneves1}.

With this proposition we can prove Lemmas \ref{increasing width} and \ref{ellipsoid width increasing}.

\begin{proof}[Proof of Lemma \ref{increasing width}]
Suppose $\width_4(S^2)<2\pi(1+\frac{1}{6})$ and choose an ellipsoid $E^2$ sufficiently close to $S^2$ so that Proposition \ref{genericellipsoid} holds.
In particular there are only $3$ almost minimising geodesic networks with length less than $2\pi(1+\frac{1}{4})$ in $E^2$ (namely, the three principal geodesics).
Indeed, any such geodesic network must have density less than $2$ by Theorem \ref{densityellipsoid} and the almost minimising condition excludes triple junctions.
Note also that $\width_4(E^2)<2\pi(1+\frac{1}{3})$.

We claim that there exist an optimal sweepout for $\width_4(E^2)$.
If that is not the case we would be able to produce a sequence of sweepouts $F_i$ with no concentration of mass such that $L[\Pi_{F_{i+1}}]<L[\Pi_{F_i}]<2\pi(1+\frac{1}{3})$. 
Thus, each $F_i$ would give us a distinct almost-minimising geodesic network with length less than $2\pi(1+\frac{1}{3})$ (Theorem \ref{existence almost minimising}), which is a contradiction.

It follows that there exists an almost-minimising geodesic network $V$ such that $\|V\|(E^2)=\width_4(E^2)$.
Thus, $V$ must be one of the three principal geodesics which implies that $\width_k(E^2)=\width_{k+1}(E^2)$ for some $k=1,2,3$.
This is a contradiction because Proposition \ref{equality width} would imply the existence of infinitely many almost-minimising geodesic networks with length $\width_k(E^2)$ and we already know that this is not possible.

We conclude that our initial assumption is false, thus $\width_4(S^2)>2\pi$.
\end{proof}

The next proof is very similar to the previous one.

\begin{proof}[Proof of Lemma \ref{ellipsoid width increasing}]
If the ellipsoid is sufficiently close to $S^2$ then we can assume that $\width_i(E^2)<2\pi(2+\frac{1}{4})$, by Theorem \ref{width sphere}.
As we have already remarked, there are only $9$ almost-minimising geodesic networks with mass less than $2\pi(2+\frac{1}{4})$ (namely the $W_j$ previously described).
If we had equality $\width_{i}=\width_{i+1}$ for any $i=1,\ldots,7$ then Proposition \ref{equality width} would give us infinitely many almost-minimising geodesic networks with mass $\width_i$, which is a contradiction.
\end{proof}

\end{appendices}